\documentclass[a4paper,10pt]{article}
\usepackage{amssymb}% \mathbb, \mathrm, \mathfrak
\usepackage{latexsym}% \leq, \geq

\begin{document}

\begin{center}
\LARGE { Moduli Spaces of Instantons on Noncommutative 4-Manifolds } \\
\vspace{2cm}
\Large { Hiroshi TAKAI} \\
\vspace{5mm}
\Large { Department of Mathematics, \\
Tokyo Metropolitan University} \\
\end{center}

\vspace{3cm}

\begin{center}
\Large  Abstract 
\end{center}

\large  \quad Studied are moduli spaces of self dual or anti-self dual 
connections on noncommutative 4-manifolds, especially deformation quantization 
of compact spin Riemannian 4-manifolds and their isometry groups have 2-torus 
subgroup. Then such moduli spaces of irreducible modules associated with highestweights of compact connected semisimple Lie groups are smooth manifolds with dimension determined by their weights. This is a generalization of 
Atiyah-Hitchin-Singer's classical result as well as Landi-Suijlekom's noncommutative 4-sphere case.

\newpage
\Large {\bf{\S1.~Introduction}}\large \quad Among many important topics of 
Yang-Mills theory, Atiyah-Hitchin-Singer [AHS] showed that the moduli spaces of  irreducible self-dual connections of principal buldles over compact spin 
4-manifolds are smooth manifolds with their dimension determined in terms of 
topological invariants related to ambient bundles. Recent development of 
noncommutative geometry serves many powerful devices toward determining a 
noncommutative analog of their result cited above. For example, Connes-Landi[CL] showed an existence of an isospectral deformation of compact spin Riemannian 
manifolds whose isometry groups have at least 2-torus subgroup. Related with 
this results, Connes-Violette [CV] observed that a deformation 
quantization of manifolds could be viewed as the fixed point algebra of the set of all smooth functions from manifolds to noncommutative torus under certain 
action of torus. Suceedingly, Landi-Suijlekom [LS] computed concerning 
the noncommutative 4-sphere, the index of the Dirac operator on certain spin 
modules as well as the dimension of the instantons of an irreducible module. 
   In this paper, we show that given a principal bundle over compact Riemannian spin 4-manifold whose isometry group has 2-torus subgroup, the moduli spaces of  instantons (anti-instantons) of irreducible projective modules associated with  highest weights are smooth manifolds with their dimension determined by ambient principal bundles and highest weights. \\
   
\Large{\bf{\S2.~Noncommutative Yang-Mills Theory}} \\
\large In this section, 
we review briefly both commutative Yang-Mills theory due to Atiyah-Hitchin-
Singer [AHS] and the noncommutative Yang-Mills theory due to Landi-Suijlekom 
[LS]. Let $M$ be a  compact spin Riemannian 4-manifold, $G$ a compact connected semisimple Lie group, and $P$ a principal $G$-bundle over $M$. 
Suppose there exist a smooth action $\varphi$ 
from 2-torus $T^2$ into the isometry group $Iso(M,g)$ of $M$ with respect to its Riemannian metric $g$. By Connes-Violette [CV], let $M_{\theta}$ be a deformation quantization of $M$ along $\theta$ as smooth *-algebra. 
Actually, there exists a quantization map 
$\mathrm{L}_{\theta}$ from $C^{\infty}(M)$ onto $M_{\theta}$. 
Then it is identified with the fixed point algebra of the set of all smooth maps from M to the noncommutative 2-torus $T^2_{\theta}$ 
under the diagonal tensor action $\varphi \otimes \alpha^{-1}$ of $T^2$ 
where $\alpha$ is the gauge action of $T^2$ on $T^2_{\theta}$. 
Since $M$ is  4-dimensional, its Hodge *-operation on the Grassmann algebra 
$\Omega(M)$ of all forms of $M$ could be shifted on the all forms 
$\Omega(M_{\theta})$ of $M_{\theta}$, 
which is denoted by $*_{\theta}$. Let us take a noncommutative principal bundle  as follows: let 
\[ G \stackrel{\varrho} \rightarrow N \stackrel{\pi} \rightarrow M \]
\noindent
be a principal $G$-bundle over $M$. Suppose there exists a smooth action 
$\tilde{\varphi}$ of a covering group $\tilde{T^2}$ of $T^2$ to 
$\mathrm{Iso}(N,\pi^*(g))$ commuting with $\varrho$, then it follows from [LS] 
that there exists a smooth action $\tilde{\varrho}$ of $G$ on 
$N_{\tilde{\theta}}$ such that
\[G \stackrel{\tilde{\varrho}} \rightarrow N_{\tilde{\theta}}
\stackrel{\pi_{\theta}} \rightarrow M_{\theta}  \]
\noindent
is a noncommutative $G$-bundle over $M_{\theta}$ where $\pi^*(g)$ means the 
pull back Riemannian metric of $g$ on $N$ under $\pi$. Let $\sigma$ be a highest weight of $G$ with respect to its maximal torus $T$ and $V_{\sigma}$ its 
irreducible $G$-module. Let $\Xi_{\sigma}$ be the fixed point algebra of 
$N_{\tilde{\theta}}\otimes V_{\sigma}$ under the diagonal tensor action 
$\tilde{\varrho}\otimes \sigma^{-1}$ of $G$. Then it is a finitely generated 
irreducible right $M_{\theta}$-module, and there exists a natural number $n$ 
and a projection $P_{\sigma}$ in $M_n(M_{\theta})$ such that 
$\Xi_{\sigma}=P_{\sigma}(M_{\theta}^n)$. We now take the Grassman connection 
$\nabla_{\sigma}=P_{\sigma}d_{\theta}$ of $\Xi_{\sigma}$ where $d_{\theta}$ is 
the canonical outer derivative of $\Omega(M_{\theta})$ derived by the ordinary 
outer derivative $d$ of $\Omega(M)$. We have nothing to know at the moment 
whether $\nabla_{\sigma}$ is self-dual or anti selfdual although this is the 
case in the special set up (for instance [CDV],[LS]). Let us denote by 
$\mathcal{C}_{\pm}(\Xi_{\sigma})$ the set of all self-dual(+), anti self-dual(-) compatible connections of $\Xi_{\sigma}$ respectively. In what follows, 
we treat the case that $\mathcal{C}_{\pm}(\Xi_{\sigma})$ are non empty although we have  nothing to check their existence for a given principal $G$-bundle over $M$ and an irreducible representation of $G$. According to [LS], the set 
$\mathcal{C}_{\mathrm{YM}}(\Xi_{\sigma})$ of all Yang-Mills connections of 
$\Xi_{\sigma}$ contains the sum of $\mathcal{C}_{+}(\Xi_{\sigma})$ and 
$\mathcal{C}_{-}(\Xi_{\sigma})$.\\
 
\Large{\bf{\S3.~Geometry of Noncommutative Instantons}} \large In this section, we analize a geometric structure of $\mathcal{C}_{\pm}(\Xi_{\sigma})$ defined in the previous section under the case where they are non empty. Let $\nabla_{\sigma} \in \mathcal{C}_{\pm}(\Xi_{\sigma})$ for an irreducible representation 
$\sigma$ of $G$ on a finite dimensional $\mathbb{R}$-vector space $V_{\sigma}$. Let $P_{\pm}$ be the projection from 
$\Omega^2(M_{\theta},\hat{\Xi}_{\sigma})$ onto 
$\Omega^2_{\pm}(M_{\theta},\hat{\Xi}_{\sigma})$, where 
$\hat{\Xi}_{\sigma}=\mathrm{End}_{M_{\theta}}(\Xi_{\sigma})$. 
Since $\mathrm{End}_{M_{\theta}}(\Xi_{\sigma})$ is isomorphic to 
$(N_{\tilde{\theta}}\otimes \mathcal{L}(V_{\sigma}))^{\tilde{\varrho}
\otimes Ad(\sigma)(G)}$, 
then its skew-adjoint part $\hat{\Xi}^{sk}_{\sigma}$ is identified with 
$(N^{s}_{\tilde{\theta}}\otimes \mathcal{L}^{sk}(V_{\sigma}))^{\tilde{\varrho}
\otimes Ad(\sigma)(G)}$, 
where $\mathcal{L}(V_{\sigma})$ is the set of all $\mathbb{C}$-linear maps 
on $V_{\sigma}$,~$\mathcal{L}^{sk}(V_{\sigma})$ its skew-adjoint part,~and 
$N^{s}_{\theta}$ is the set of all self-adjoint elements of 
$N_{\tilde{\theta}}$.
Since $\mathcal{L}(V_{\sigma})=M_{n_{\sigma}}(\mathbb{C})$, 
then it follows that
$\{T\in \mathcal{L}^{sk}(V_{\sigma})|Tr(T)=0 \}=\mathcal{SU}(n_{\sigma})$ 
where $n_{\sigma}=\mathrm{dim}_{\mathbb{C}}V_{\sigma}$ and $\mathcal{SU}(n_{\sigma})$ is the Lie algebra of $\mathrm{SU}(n_{\sigma})$. We then deduce that 
\[M_{n_{\sigma}}(\mathbb{C})=\mathbb{C}\mathrm{I}_{n_{\sigma}} 
\oplus \mathcal{SU}(n_{\sigma})_{\mathbb{C}} \]
\noindent
where $\mathcal{SU}(n_{\sigma})_{\mathbb{C}}$ is the complexification of $\mathcal{SU}(n_{\sigma})$. Therefore it implies that 
\[\hat{\Xi}_{\sigma}=M_{\theta}\oplus_{M_{\theta}} (N_{\tilde{\theta}} 
\otimes \mathcal{SU}(n_{\sigma})_{\mathbb{C}})^{\tilde{\varrho}\otimes Ad(\sigma)(G)} \]
\noindent
Let us define
\[\Gamma(ad_{\sigma}(N_{\tilde{\theta}}))=(N_{\tilde{\theta}}\otimes \mathcal{SU}(n_{\sigma})_{\mathbb{C}})^{\tilde{\varrho}\otimes Ad(\sigma)(G)}. \]
\noindent
Moreover, we put
\[\Omega^0(ad_{\sigma}(N_{\tilde{\theta}}))
=\Gamma(ad_{\sigma}(N_{\tilde{\theta}})),
\Omega^1(ad_{\sigma}(N_{\tilde{\theta}}))
=\Omega^1(M_{\theta},\Gamma(ad_{\sigma}(N_{\tilde{\theta}}))) \]
\noindent
,~and 
\[\Omega^2_{\pm}(ad_{\sigma}(N_{\tilde{\theta}}))=
\mathrm{P}_{\pm}\Omega^2(M_{\theta},\Gamma(ad_{\sigma}(N_{{\tilde{\theta}}})))
\] 
\noindent
respectively. We then introduce an inner product $<~|~>$ on 
$\Omega^j(ad_{\sigma}(N_{\tilde{\theta}}))~(j=0,1,2)$ by 
\[ <~\omega~|~\eta~>=\int_{M_{\theta}}\mathrm{tr}(*_{\theta}(\omega^* *_{\theta}\eta)) \]
\noindent
where $*_{\theta}$ is the Hodge operation on $\Omega^*(M_{\theta})$,and $\mathrm{tr}$ is the canonical trace on $M_{n_{\sigma}}(\mathbb{C})$. Using this, 
we also induce the metric topology on $\mathcal{C}_{\pm}(\Xi_{\sigma})$ 
by the following lemma:\\

\large{\bf{Lemma 3.1}}~~Let $\nabla_{\sigma} \in \mathcal{C}_{\pm}(\Xi_{\sigma})$. Then it follows that 
\[\mathcal{C}_{\pm}(\Xi_{\sigma})
=\nabla_{\sigma}+\Omega^1_{\pm}(ad_{\sigma}(N_{\tilde{\theta}}) \]
\noindent
where ~$\Omega^1_{\pm}(ad_{\sigma}(N_{\tilde{\theta}}))$ is 
the set of all $\omega \in \Omega^1(ad_{\sigma}(N_{\tilde{\theta}}))$ 
satisfying the equation:
\[  \hat{\nabla}^{\mp}(\omega)+(\omega^2)_{\mp}=0 ~,\]
\noindent
and $(\omega^2)_{\mp}=\mathrm{P}_{\mp}(\omega^2)$ respectively.\\

Proof. Let $\nabla \in \mathcal{C}_{\pm}(\Xi_{\sigma})$ and put $\omega=\nabla-\nabla_{\sigma}$, then it implies by definition that 
$\omega \in \hat{\Xi}_{\sigma}$. Since $\nabla,~\nabla_{\sigma} \in \mathcal{C}_{\pm}(\Xi_{\sigma})$, it follows that
\[\mathrm{P}_{\mp}[\nabla,[\nabla,X]]=\mathrm{P}_{\mp}[F_{\nabla},X]=[F^{\mp}_{\nabla},X]=0 \] 
\noindent
for all $X \in \hat{\Xi}_{\sigma}$ respectively, and the similar statement holds for $\nabla_{\sigma}$. We compute that 
\[F_{\nabla}=F_{\nabla_{\sigma}}+[\nabla_{\sigma},\omega]+\omega^2 \], 
\noindent
which implies the conclusion. ~~~ Q.E.D.  \\

As $M$ is a spin manifold, it follows from [CV] that there exists a 
noncommutative spin structure $S_{\theta}$ of $M_{\theta}$ 
as a right $M_{\theta}$-module induced by the given spin structure $S$ of $M$. 
Let $S^{\pm}_{\theta}$ be the half spin structures of $M_{\theta}$ induced by 
the given half spin structures  $S^{\pm}$ of $M$ respectively. We use the same 
notations as their smooth sections on $M_{\theta}$. 
Then the following statement is easily seen by their definition: \\

\large{\bf{Lemma 3.2}}~~
\[(1)~~\Omega^1(M_{\theta})=S^{\pm}_{\theta}\otimes S^{\mp}_{\theta}\]
\noindent
\[~~~~~~~~~(2)~~M_{\theta}\oplus \Omega^2_{\pm}(M_{\theta})=
S^{\pm}_{\theta} \otimes S^{\pm}_{\theta} \]
\noindent
respectively.\\
 
Proof.~(1): ~By definition, $S^{\pm}_{\theta}$ is the algebra of all 
even (odd) polynomials of $\Omega^1(M_{\theta})$ with the Clifford 
multiplication $c_{\theta}$ respectively. 
Then $c_{\theta}(\Omega^1(M_{\theta})S^{\pm}_{\theta}=S^{\mp}_{\theta}$ 
respectively, which implies the conclusion. ~(2):~The result is true if 
$\theta=0$. This isomorphism can also be chosen $\tilde{T}^2$-equivariantly, 
where $\tilde{T}^2$ is a covering group of $T^2$. Then the statement follows. 
~~~Q.E.D.\\

\noindent
By Lemmas 3.1, $\mathcal{C}_{\pm}(\Xi_{\sigma})$ is identified with 
$\Omega^1_{\pm}(ad_{\sigma}(N_{\tilde{\theta}}))$ as an affine space 
respectively. We then show the next three lemma which seems to be quite useful 
showing our main result:\\

\large{\bf{Lemma 3.3}}~~Let $T_{\nabla}(\mathcal{C}_{\pm}(\Xi_{\sigma}))$ be 
the tangent space of $\mathcal{C}_{\pm}(\Xi_{\sigma})$ at 
$\nabla \in \mathcal{C}_{\pm}(\Xi_{\sigma})$. Then it follows that 
\[T_{\nabla}(\mathcal{C}_{\pm}(\Xi_{\sigma}))=\{\omega \in 
\Omega^1(ad_{\sigma}(N_{\tilde{\theta}}))~|~\mathrm{P}_{\mp}[\nabla,\omega]=0 \}~,\] 
\noindent
respectively. \\

Proof.~~Let $\omega_t=\nabla_t-\nabla$ for a smooth curve $\nabla_t \in 
\mathcal{C}_{\pm}(\Xi_{\sigma})$ with $\nabla_0=\nabla$. 
Then it follows from definition that 
$\omega_t \in \Omega^1(ad_{\sigma}(N_{\tilde{\theta}}))$. We put $\omega$ 
the derivative of $\omega_t$ at $t=0$ in 
$\Omega^1(ad_{\sigma}(N_{\tilde{\theta}}))$. 
Then we see that
\[F_{\nabla_t}=F_{\nabla}+[\nabla,\omega_t]+\omega_t^2 ~,\] 
\noindent
which deduces $F_{\nabla_t}\prime(0)=[\nabla,\omega]$ 
taking their derivatives at $t=0$ since $\omega_0=0$.
As $\nabla_t \in \mathcal{C}_{\pm}(\Xi_{\sigma})$, it implies that 
$F_{\nabla_t} \in \Omega^2_{\pm}(ad(N_{\tilde{\theta}}))$, which means that 
$\mathrm{P}_{\mp}[\nabla,\omega]=0$ respectively. ~~~Q.E.D. \\

\noindent
Let $\Gamma(\Xi_{\sigma})$ be the gauge group acting on the set 
$\mathcal{C}(\Xi_{\sigma})$ of all compatible connections of 
$\Xi_{\sigma}$ by 
$\gamma_u(\nabla)=u{\nabla}u^*,(u \in \Gamma(\Xi_{\sigma},\nabla \in \mathcal{C}(\Xi_{\sigma}))$,~and 
$\Gamma(\Xi_{\sigma}) \cdot \nabla$ the orbit of $\nabla$ under $\Gamma(\Xi_{\sigma})$. Then we show the following lemma:\\

\large{\bf{Lemma 3.4}}~~~~~Let 
$\mathrm{T}_{\nabla}(\Gamma(\Xi_{\sigma})\cdot\nabla)$ be the tangent space 
of $\Gamma(\Xi_{\sigma})\cdot\nabla$ at 
$\nabla \in \mathcal{C}(\Xi_{\sigma})$. 
Then we have that
\[\mathrm{T}_{\nabla}(\Gamma(\Xi_{\sigma})\cdot\nabla)=\{~[\nabla,X]~|~X \in 
\Omega^0(ad(N_{\tilde{\theta}}))~\} \]
\noindent
for any $\nabla \in \mathcal{C}_{\pm}(\Xi_{\sigma})$ respectively.\\

Proof.~~Let $X \in \Omega^0(ad_{\sigma}(N_{\tilde{\theta}}))$ and 
put $\varphi_t=e^{-tX} \in \Gamma(\Xi_{\sigma})$ 
for all $t \in \mathbb{R}$. Then we know that 
\[\gamma_{\varphi_t}(\nabla)(\xi)=e^{-tX}\cdot\nabla(e^{tX}\xi)=t[\nabla,X](\xi)+\nabla(\xi) \]
\noindent
for all $\xi \in \Xi_{\sigma}$. Taking their derivatives at $t=0$, we have that \[ \frac{d}{dt}\gamma_{\varphi_t}(\nabla)\Big{|}_{t=0}(\xi)=[\nabla,X](\xi)\]
\noindent
for all $\xi \in \Xi_{\sigma}$. Therefore the conclusion follows.~~~Q.E.D. \\

\noindent
As we know that 
$\Gamma(\Xi_{\sigma})\cdot\nabla \subseteq \mathcal{C}_{\pm}(\Xi_{\sigma})$ 
for any $\nabla \in \mathcal{C}_{\pm}(\Xi_{\sigma})$ respectively, 
we can imagine by Lemma 3,3 and 3.4 the next corollary: \\

\large{\bf{Corollary 3.5}}~~Let $\mathcal{M}_{\pm}(\Xi_{\sigma})$ be the moduli  space of $\mathcal{C}_{\pm}(\Xi_{\sigma})$ by $\Gamma(\Xi_{\sigma})$. Then the  tangent space$T_{[\nabla]}(\mathcal{M}_{\pm}(\Xi_{\sigma}))$ of 
$\mathcal{M}_{\pm}(\Xi_{\sigma})$ at 
$[\nabla] \in \mathcal{M}_{\pm}(\Xi_{\sigma})$ is isomorphic to 
$T_{\nabla}(\mathcal{C}_{\pm}(\Xi_{\sigma}))/\mathrm{T}_{\nabla}(\Gamma(\Xi_{\sigma})\cdot\nabla)$ 
as a $\mathbb{C}$-linear space respectively.\\

\noindent
Let us define $\hat{\nabla}^{\pm}(\omega)=\mathrm{P}_{\pm}[\nabla,\omega]$ 
for all $\omega \in \Omega^1(ad_{\sigma}(N_{\tilde{\theta}}))$ respectively. 
We then can show the following lemma which is well known in undeformed cases 
by [AHS] and in a special deformed case by [LS]:\\

\large{\bf{Lemma 3.6}}~~
\[0 \rightarrow \Omega^0(ad_{\sigma}(N_{\tilde{\theta}}))
 \stackrel{\hat{\nabla}} 
\rightarrow \Omega^1(ad_{\sigma}(N_{\tilde{\theta}})) 
\stackrel{\hat{\nabla}^{\pm}} \rightarrow 
\Omega^2_{\pm}(ad_{\sigma}(N_{\tilde{\theta}})) \rightarrow 0 \]
\noindent
for $\nabla \in \mathcal{C}_{\mp}(\Xi_{\sigma})$, where $
\hat{\nabla}(X)=[\nabla,X]$ for $X \in 
\Omega^0(ad_{\sigma}(N_{\tilde{\theta}}))$ and 
$\hat{\nabla}^{\pm}(\omega)=\mathrm{P}_{\pm}[\nabla,\omega]$ for 
$\omega \in \Omega^1(ad_{\sigma}(N_{\tilde{\theta}}))$ respectively. \\

Proof.~~We first show that 
$\hat{\nabla}^{\pm}\cdot \hat{\nabla}=0$. Indeed, 
as $\nabla \in \mathcal{C}_{\mp}(\Xi_{\sigma})$,

\[\hat{\nabla}^{\pm}\cdot \hat{\nabla}(X)=
P_{\pm}[\nabla,[\nabla,X]]=P_{\pm}[F_{\nabla},X]=
[F^{\pm}_{\nabla},X]=0\]
\noindent
for any $X \in \Omega^0(ad_{\sigma}(N_{\tilde{\theta}}))$ respectively. 
We next show that$\mathrm{Ker}~\hat{\nabla}=0$. In fact, 
suppose $X \in \mathrm{Ker}~\hat{\nabla}$. Since $\hat{\nabla}$ commutes with 
the action $\alpha$ of $T^2$ and $ad(N_{\tilde{\theta}})$ has 
$\mathcal{SU}(n_{\sigma})$ as fibres, then it follows that 
\[\hat{\nabla}=d_{\tilde{\theta}}+\omega_n \] 
\noindent
for some $\omega_n \in M_n(\Omega^1(N_{\tilde{\theta}}))$, 
where $(\omega_n)_{j,k}$ are all central in $\Omega^1(N_{\tilde{\theta}})$, 
which means that $\mathrm{L}_{\theta}(\omega_n)=\omega_n$. Hence it implies that $\hat{\nabla}(X)=\mathrm{P}_{\sigma}\cdot d_{\theta}(X)=0$. 
Since $X=\mathrm{L}_{\theta}(X^0)$ for a $X^0 \in \Omega^0(ad(N))$, 
then we have that 
\[ \hat{\nabla}(X)=\mathrm{P}_{\sigma}\cdot d_{\theta}(X)=\mathrm{L}_{\theta}
\cdot \mathrm{p}_{\sigma}\cdot d_M(X^0)=0~, \]
\noindent
where $\mathrm{p}_{\sigma}$ is the projection in $M_n(C^{\infty}(M))$ 
corresponding to the right $C^{\infty}(M)$-module $\Xi=
(C^{\infty}(N)\otimes V_{\sigma})^{\varrho \otimes \sigma(G)}$ and $d_M$ is 
the outer derivative of $M$. Then $\mathrm{p}_{\sigma}\cdot d_M(X^0)=0$, 
which implies that $X^0=0 \in \Omega^0(ad(N))$ by the irreducibility of 
$\mathrm{p}_{\sigma}\cdot d_M$ on $\Xi$. This means that 
$\mathrm{Ker}~\hat{\nabla}=0$.  ~~~~~Q.E.D.  \\

\noindent
Let us define ~$\mathbb{H}^1_{\pm}
=\mathrm{Ker}~\hat{\nabla}^{\pm}/\mathrm{Im}~\hat{\nabla}$~and~
$\mathbb{H}^2_{\pm}
=\Omega^2_{\pm}(ad_{\sigma}(N_{\tilde{\theta}})/\mathrm{Im}~\hat{\nabla}^{\pm}$. It follows from Lemma 3.5 that $\mathbb{H}^1_{\pm}$ is isomorphic to 
$T_{[\nabla]}(\mathcal{M}_{\pm}(\Xi_{\sigma}))$ as a $\mathbb{C}$-linear space 
for any $\nabla \in \mathcal{C}_{\mp}(\Xi_{\sigma})$ respectively. Moreover, we  introduce the Laplace type $M_{\theta}$-operators $\Delta^j_{\pm}$ on 
$\Omega^j(ad_{\sigma}(N_{\tilde{\theta}}))$ \\
\noindent
for $j=1,2$ as follows: 
let $\Delta^{\pm}_j=(\hat{\nabla}^{\pm})^*\cdot \hat{\nabla}^{\pm}+\hat{\nabla}\cdot (\hat{\nabla})^*$ on $\Omega^1(ad_{\sigma}(N_{\tilde{\theta}}))$ and 
$\Delta^{\pm}_j=\hat{\nabla}^{\pm}\cdot (\hat{\nabla}^{\pm})^*$ on 
$\Omega^2(ad_{\sigma}(N_{\tilde{\theta}}))$ respectively. 
Then we easily observe the following lemma: \\

\large{\bf{Lemma 3.7}}~~$\mathbb{H}^j_{\pm}$ are isomorphic to 
$\mathrm{Ker}~\Delta^{\pm}_j~(j=1,2)$ 
as a right $M_{\theta}$-module respectively, and they are finitely generated 
projective right $M_{\theta}$-modules respectively. \\

Proof.~~$\mathbb{H}^1_{\pm}$ is $M_{\theta}$-isomorphic to 
$\mathrm{Ker}~\hat{\nabla}^{\pm} \cap \mathrm{Ker}~(\hat{\nabla})^*$, which is 
equal to $\mathrm{Ker}~\Delta^{\pm}_j$. The similar way is also valid for 
$\mathbb{H}^2_{\pm}$. As $\Delta^{\pm}_j$ are elliptic, the rest is well known.
~~~Q.E.D. \\

\noindent
By Lemma 3.6, the elliptic complex can be described by the following generalized signature $M_{\theta}$-operator:
\[ \hat{\nabla}^{\pm}+(\hat{\nabla})^*:\Omega^1(ad_{\sigma}(N_{\tilde{\theta}})) \rightarrow \Omega^0(ad_{\sigma}(N_{\tilde{\theta}})) 
\oplus \Omega^2_{\pm}(ad_{\sigma}(N_{\tilde{\theta}})) \]
\noindent
for all $\nabla \in \mathcal{C}_{\mp}(\Xi_{\sigma})$. By Lemma 3.7, $\mathbb{H}^j_{\pm}$ induce the $\mathbb{K}_0(M_{\theta})$-element $[\mathbb{H}^j_{\pm}]~(j=1,2)$ whose ranks are finite. Let us define the K-theoretic index of $\hat{\nabla}^{\pm}+(\hat{\nabla})^*$ as follows:
\[\mathrm{Index}_{M_{\theta}}(\hat{\nabla}^{\pm}+(\hat{\nabla})^*)=[\mathrm{Ker}(\hat{\nabla}^{\pm}+(\hat{\nabla})^*)]-[\mathrm{Coker}(\hat{\nabla}^{\pm}+(\hat{\nabla})^*)] \]
\noindent
Using Lemma 3.6, we then obtain the following lemma which is a noncommutative 
version of Atiyah-Singer index theorem due to Connes [C]~(cf:[MS]): \\

\large{\bf{Lemma 3.8}}~~
\[ \mathrm{Index}_{M_{\theta}}(\hat{\nabla}^{\pm}+(\hat{\nabla})^*)
=\mp~[\mathbb{H}^1_{\pm}]~\pm[\mathbb{H}^2_{\pm}] \]
\noindent
respectively. \\

\noindent
Let $[M_{\theta}]$ be the fundamental cyclic cocycle of $M_{\theta}$ which is 
essentially defined by the JLO-cocycle appeared in [C]. 
Then it follows by Lemma 3.8 that \\

\large{\bf{Corollary 3.9}}~~
\[<[M_{\theta}],\mathrm{Index}_{M_{\theta}}(\hat{\nabla}^{\pm}+(\hat{\nabla})^*)>
=\mp\mathrm{rank}_{M_{\theta}}[\mathbb{H}^1_{\pm}]\pm\mathrm{rank}_{M_{\theta}}[\mathbb{H}^2_{\pm}] \]
\noindent
respectively. \\

\noindent
On the other hand, we consider the Dirac type $M_{\theta}$-operator 
$D^{\pm}_{\nabla}$ on $\Omega^1(ad_{\sigma}(N_{\tilde{\theta}}))$ 
associated with $\nabla \in \mathcal{C}_{\pm}(\Xi_{\sigma})$ 
defined by the following process: By Lemma 3.2.(1), \\

$\Omega^1(ad_{\sigma}(N_{\tilde{\theta}})
=\Omega^0(ad_{\sigma}(N_{\tilde{\theta}}))
\otimes (S^{\pm}_{\theta}\otimes S^{\mp}_{\theta})$ \\

$\stackrel{\hat{\nabla}} \rightarrow 
\Omega^0(ad_{\sigma}(N_{\tilde{\theta}}))\otimes \Omega^1(M_{\theta})
\otimes (S^{\pm}_{\theta}\otimes S^{\mp}_{\theta}) $ \\

$\stackrel{\mathrm{Id} \otimes c^{\pm}_{\theta}} \rightarrow 
\Omega^0(ad_{\sigma}(N_{\tilde{\theta}}))\otimes 
(S^{\mp}_{\theta}\otimes S^{\mp}_{\theta})=\Omega^0(ad_{\sigma}(N_{\tilde{\theta}})) \oplus \Omega^2_{\mp}(ad_{\sigma}(N_{\tilde{\theta}}))$ \\

\large{\bf{Lemma 3.10}}~~Given a $\nabla \in \mathcal{C}_{\pm}(\Xi_{\sigma})$,
\[{\mathrm{Index}}_{M_{\theta}}(D^{\pm}_{\nabla})={\mathrm{Index}}_{M_{\theta}}(\hat{\nabla}^*+\hat{\nabla}^{\mp}) \] 
\noindent 
in $K_0(M_{\theta})$ respectively. \\ 

Proof.~~Since $\Omega(M_{\theta})\simeq \Omega(M)$ as vector spaces and 
$\hat{\nabla}$ and $*$-operation commute with $T^2$-action $\alpha$, then 
the operator $\hat{\nabla}^*+\hat{\nabla}^{\mp}$ can be shifted to the elliptic  operator $(d_M)^*+p_{\mp}\cdot d_M$ from $\Omega^1(M)$ to 
$\Omega^0(M) \oplus \Omega^2_{\mp}(M)$. Let $P^0_{\mp}$ be the projection from 
$\Omega^2(M)$ to $\Omega^2_{\mp}(M)$ respectively. Using the same argument as 
in [AHS] and taking their principal symbols, it follows from Lemma 3.2 that 
$(d_M)^*+P^0_{\mp}\cdot d_M$ is identified with the Dirac operator $D^{\pm}$ 
in the following:
\[D^{\pm}:~S^{\pm} \otimes S^{\mp} \rightarrow S^{\mp} \otimes S^{\mp} \]~,
\noindent
so far taking their $C^{\infty}(M)$-indices. Then 
$\hat{\nabla}^*+\hat{\nabla}^{\mp}$ is identified with $D^{\pm}_{\sigma}$ 
defined before:
\[ D^{\pm}_{\sigma}: \Omega^0(ad_{\sigma}(N_{\tilde{\theta}}))
\otimes (S^{\pm}_{\theta}\otimes S^{\mp}_{\theta}) \rightarrow 
\Omega^0(ad_{\sigma}(N_{\tilde{\theta}}))
\otimes (S^{\mp}_{\theta}\otimes S^{\mp}_{\theta}) \]
\noindent
so far taking their $M_{\theta}$-indices. This implies the conclusion.~Q.E.D.\\

\noindent 
In what follows, we want to determine the geometric structure of the moduli 
space $\mathrm{M}_{\pm}(\Xi_{\sigma})$ of $\mathrm{C}_{\pm}(\Xi_{\sigma})$ by 
the gauge group $\Gamma(\Xi_{\sigma})$. First of all, we introduce 
a noncommutative Kuranishi map in the following process : 
Given a $\nabla \in \mathcal{C}_{\mp}(\Xi_{\sigma})$, 
we define a Sobolev norm $||\cdot||_p,~(p\gg 4)$ on 
$\Omega^1_{\pm}(ad_{\sigma}(N_{\tilde{\theta}})$ by \\

$~~~~~~~~~~~~~~~~~~||\omega||_p = \sum^p_{k=0}\int_{M_{\theta}}~
||(\hat{\nabla})^k(\omega)||^2$ \\

\noindent
for $\omega \in \Omega^1_{\pm}(ad_{\sigma}(N_{\tilde{\theta}}))$. 
Let $S^p_{\pm}$  be the completion of 
$\Omega^1_{\pm}(ad_{\sigma}(N_{\tilde{\theta}}))$ with respect to $||\cdot||_p$. By their definition, $\hat{\nabla},~\hat{\nabla}^*$ are bounded from 
$S^p_{\pm}$ to $S^{p+1}_{\pm},~S^{p-1}_{\pm}$ respectively. \\

\noindent 
Since  $\Delta^{\pm}_j$ is elliptic in $\mathrm{End}_{M_{\theta}}(S^p_{\pm})$, 
it is a $M_{\theta}$-Fredholm operator on $S^p_{\pm}$ as well. 
 Let $\mathrm{P}_{\mathbb{H}^2_{\pm}}$ be the projection on $\mathbb{H}^2_{\pm}$  ,and define the Green $M_{\theta}$-operator $G^{\pm}$ of $\Delta^{\pm}_2$ 
 on $S^p_{\pm}$ as follows :
\[G^{\pm}=(\mathrm{Id}_{S^p_{\pm}}-\mathrm{P}_{\mathbb{H}^2_{\pm}})\cdot 
(\Delta^{\pm}_2)^{-1}= (\Delta^{\pm}_2)^{-1} \cdot 
(\mathrm{Id}_{S^p_{\pm}}-\mathrm{P}_{\mathbb{H}^2_{\pm}}) \]
\noindent
respectively. Then we easily see that \\

$(1):~\mathrm{Im}~G^{\pm}=\mathrm{Im}~\Delta^{\pm}_2$, \\

$(2):~G^{\pm} \cdot \Delta^{\pm}_2=\Delta^{\pm}_2 \cdot G^{\pm}$, \\

$(3):~G^{\pm} \cdot \Delta^{\pm}_2=
\mathrm{Id}_{S^p_{\pm}}-\mathrm{P}_{\mathbb{H}^2_{\pm}}$ \\

\noindent
We now introduce a densely defined map $\Pi^{\pm}_{\nabla}$ 
on $S^p_{\pm}$ by
\[\Pi^{\pm}_{\nabla}(\omega)= 
\omega+(\hat{\nabla}^{\pm})^* \cdot G^{\pm}\{(\omega^2)_{\pm}\} \]
\noindent
for all $\omega \in \Omega^1_{\mp}(ad_{\sigma}(N_{\tilde{\theta}}))$ 
respectively, which is called a noncommutative Kuranishi map, which is known 
in commutative cases (cf:[AHS]). Then we deduce the following observation: \\

\large{\bf{Lemma 3.11}}~~\\

$~~~~~~\Omega^1_{\pm}(ad_{\sigma}(N_{\tilde{\theta}}))=
\mathrm{Ker}~(\hat{\nabla}^{\mp} \cdot \Pi^{\mp}_{\nabla}) \cap 
\{P_{\mathbb{H}^2_{\mp}}((\omega^2)_{\mp})=0 \}$ \\

\noindent
respectively.\\

Proof.~~By Lemma 3.1, $\omega \in \Omega^1_{\pm}(ad_{\sigma}(N_{\tilde{\theta}}))$ satisfies 
\[\hat{\nabla}^{\mp}(\omega)+(\omega^2)_{\mp}=0 ~,\]
\noindent
respectively. Then we compute \\

$\hat{\nabla}^{\mp}\cdot \Pi^{\mp}_{\nabla}(\omega)=\hat{\nabla}^{\mp}(\omega)+\Delta^{\mp}_2\cdot G^{\mp}((\omega^2)_{\mp})$ \\

$~~~~=\hat{\nabla}^{\mp}(\omega)+(\omega^2)_{\mp}-P_{\mathbb{H}^2_{\mp}}((\omega^2)_{\mp}=-P_{\mathbb{H}^2_{\mp}}((\omega^2)_{\mp} $ \\

\noindent
, which implies that 
\[\hat{\nabla}^{\mp}\cdot \Pi^{\mp}_{\nabla}(\omega)=-P_{\mathbb{H}^2_{\mp}}((\omega^2)_{\mp} ~.\]
\noindent
By the definition of $\mathbb{H}^2_{\pm}$, it follows that 
\[(\hat{\nabla}^{\mp})^* \cdot \hat{\nabla}^{\mp}\cdot \Pi^{\mp}_{\nabla}(\omega)=0 \]
\noindent
, which implies by taking their inner products that 
\[\hat{\nabla}^{\mp}\cdot \Pi^{\mp}_{\nabla}(\omega)
=P_{\mathbb{H}^2_{\mp}}((\omega^2)_{\mp}=0~.\] 
\noindent
respectively. The converse implication is also easily seen. ~Q.E.D. \\

\noindent
We then observe together with Lemma 3.1 the following corollary: \\

\large{\bf{Corollary 3.12}}~~If $\nabla \in \mathcal{C}_{\pm}(\Xi_{\sigma})$, 
then it follows that
\[\mathcal{C}_{\pm}(\Xi_{\sigma})=\nabla+\mathrm{Ker}~(\hat{\nabla}^{\mp} \cdot \Pi^{\mp}_{\nabla}) \cap \{P_{\mathbb{H}^2_{\mp}}((\omega^2)_{\mp})=0 \} \]
\noindent
respectively. \\

\noindent
It also follows from Lemma 3.11 that \\

\large{\bf{Corollary 3.13}}~~Suppose $\mathbb{H}^2_{\mp}=0$ for a $\nabla \in \mathcal{C}_{\pm}(\Xi_{\sigma})$, then 
it follows that
\[\Omega^1_{\pm}(ad_{\sigma}(N_{\tilde{\theta}}))=
\mathrm{Ker}~(\hat{\nabla}^{\mp}\cdot \Pi^{\mp}_{\nabla}) \]
\noindent
Moreover, we obtain by using Lemma 3.11 the following lemma: \\

\large{\bf{Lemma 3.14}}~~
\[\Pi^{\mp}_{\nabla}(\Omega^1_{\pm}(ad_{\sigma}(N_{\tilde{\theta}})))=\mathbb{H}^1_{\mp} \]
\noindent
respectively. \\

Proof.~~As $\hat{\nabla}^{\mp}\cdot \hat{\nabla}=0$ for all $\nabla \in \mathcal{C}_{\pm}(\Xi_{\sigma})$, then $(\hat{\nabla})^*\cdot (\hat{\nabla}^{\mp})^*=0$. Let $\eta=\Pi^{\mp}_{\nabla}(\omega)$ for any $\omega \in \Omega^1_{\pm}(ad_{\sigma}(N_{\tilde{\theta}}))$. Then we have that
\[(\hat{\nabla})^*(\eta)=(\hat{\nabla})^*(\omega)+\Delta^{\mp}_2\cdot G^{\mp}(
\omega^2)_{\mp}=0 \]~,
\noindent
respectively. By Lemma 3.11, ~~
\[ \hat{\nabla}^{\mp}(\eta)=P_{\mathbb{H}^2_{\mp}}((\omega^2)_{\mp})=0 \]~,
respectively. Since $\mathrm{Ker}~(\hat{\nabla})^*=\Omega^1(ad_{\sigma}(N_{\tilde{\theta}}))$, then the conclusion follows. ~~~Q.E.D. \\

\noindent
Using the Sobolev space $S^p_{\pm}$, we compute the Frechet differentiation 
$\frac{d\Pi^{\pm}_{\nabla}}{dt}\Big{|}_{t=0}$ of 
$\Pi^{\pm}_{\nabla}$ on $\Omega^1_{\pm}(ad_{\sigma}(N_{\tilde{\theta}}))$ in the following lemma: \\

\large{\bf{Lemma 3.15}}~~
\[\frac{d\Pi^{\mp}_{\nabla}}{dt}\Big{|}_{t=0}=
\mathrm{Id} \]
\noindent
on $\Omega^1_{\pm}(ad_{\sigma}(N_{\tilde{\theta}}))$ 
with respect to $||\cdot||_p$ respectively.\\ 

Proof.~~\\

$\frac{d\Pi^{\mp}_{\nabla}}{dt}\Big{|}_{t=0}(\omega)
=\displaystyle\lim_{t\to 0}~t^{-1}(\Pi^{\mp}_{\nabla}(t\omega)-\Pi^{\mp}_{\nabla}(0))$ \\
$~~~~~~~~~~~~~~~~~~~=\displaystyle\lim_{t\to 0}~(\omega+t(\hat{\nabla}^{\mp})^* \cdot G^{\mp}((\omega^2)_{\mp}))=\omega $ \\

\noindent
for any $\omega \in \Omega^1_{\pm}(ad_{\sigma}(N_{\tilde{\theta}}))$, 
which implies the conclusion. ~~~Q.E.D. \\

\noindent
By the inverse function theorem on Banach spaces, there exists a $\epsilon >0$ 
neighborhood $U(0,\epsilon)$ of $0 \in S^p_{\pm}$ with respect to $||\cdot||_p$
 on which $(\Pi^{\pm}_{\nabla})^{-1}$ is diffeomorphic. Using this fact, we next show the following lemma: \\

\large{\bf{Lemma 3.16}}~~~~~$U(0,\epsilon) \cap \mathbb{H}^1_{\pm}$ is 
diffeomorphic to $(\Pi^{\pm}_{\nabla})^{-1}(U(0,\epsilon))$ \\
\noindent
$~~\cap~\Omega^1_{\mp}(ad_{\sigma}(N_{\tilde{\theta}}))$ 
under $(\Pi^{\pm}_{\nabla})^{-1}$ respectively.\\

Proof.~~By Lemma 3.14, we know that 
\[\Pi^{\pm}_{\nabla}(\Omega^1_{\mp}(ad_{\sigma}(N_{\tilde{\theta}})))=\mathbb{H}^1_{\pm} \]
\noindent
respectively which implies the conclusion. ~~~Q.E.D. \\

\noindent 
Let $P\Gamma(\Xi_{\sigma})=\Gamma(\Xi_{\sigma})/Z(\Gamma(\Xi_{\sigma}))$ be 
the  projective gauge group of $\Gamma(\Xi_{\sigma})$ 
where $Z(\Gamma(\Xi_{\sigma})$ is the center of $\Gamma(\Xi_{\sigma})$. 
Since $\Xi_{\sigma}$ is an irreducible right $M_{\theta}$-module, 
$Z(\Gamma(\Xi_{\sigma})$ is actually the center $ZU(M_{\theta}))$ of 
the unitary group $U(M_{\theta})$ of $M_{\theta}$. We define the gauge action 
$\gamma$ of $P\Gamma(\Xi_{\sigma})$ on $\mathcal{C}_{\pm}(\Xi_{\sigma})$ by 
$\gamma_{[u]}(\nabla)(\xi)=u \nabla u^*(\xi)$ for any 
$[u]=uZU(M_{\theta})) \in P\Gamma(\Xi_{\sigma}),~\nabla \in \mathcal{C}_{\pm}(\Xi_{\sigma})$. 
We then easily see that: \\

\large{\bf{Lemma 3.17}}~~The action $\gamma$ is effective. \\

Proof.~~The statement is almost clear by its definition. ~~~Q.E.D. \\

\noindent
As the similar way to the commutative cases (cf:[AHS],[FU]), we then show 
the so-called slice theorem in noncommutative setting: \\

\large{\bf{Lemma 3.18}}~~Let $\epsilon$ be in Lemma 3.16. Then there exist \\
$0<\delta \leq \epsilon$,~$U(0,\delta)=\{\omega \in \Omega^1(ad_{\sigma}(N_{\tilde{\theta}}))~|~||\omega||_p <\delta \}$,~and a mapping 
\[ \psi:~\nabla+U(0,\delta) \rightarrow P\Gamma(\Xi_{\sigma}) \] 
\noindent
with the property that \\

$(1):~\hat{\nabla}^*(\gamma_{\psi(\nabla+\omega)}^*(\nabla+\omega)-\nabla)=0$ \\

(2):~~~there exists a $P\Gamma(\Xi_{\sigma})$-equivariant diffeomorphism 
\[~~~\Phi:\nabla+U(0,\delta) \rightarrow U(0,\mathrm{Id}) \subseteq 
\mathrm{Ker}~\hat{\nabla}^* \times P\Gamma(\Xi_{\sigma}).\]
\noindent
defined by 
\[\Phi(\nabla+\omega)=(\gamma_{\psi(\nabla+\omega)}^*(\nabla+\omega)-\nabla,\psi(\nabla+\omega)) \]

Proof.~~Let us consider the next equation defined by 
\[ \mathrm{L}_{\hat{\nabla}}(\omega,[u])=\hat{\nabla}^*(u^*{\nabla}u+u^*{\omega}u)=0 \]
\noindent
for $\nabla+\omega \in \mathcal{C}_{\mp}(\Xi_{\sigma})$ and $[u] \in P\Gamma(\Xi_{\sigma})$. 
Then the differential $\delta \mathrm{L}_{\hat{\nabla}}$ of 
$\mathrm{L}_{\hat{\nabla}}$ at $(0,\mathrm{Id})$ is given by 
\[\delta \mathrm{L}_{\hat{\nabla}}:\Omega^1_{\mp}(ad_{\sigma}(N_{\tilde{\theta}})) \times \Omega^0(ad_{\sigma}(N_{\tilde{\theta}})) \rightarrow \Omega^0(ad_{\sigma}(N_{\tilde{\theta}})) \]
\[\delta \omega \oplus \delta [u] \rightarrow \hat{\nabla}^*([\nabla,\delta u]+\delta \omega) \]~,
\noindent
The partial differential $\delta_2\mathrm{L}_{\hat{\nabla}}$ of $\delta \mathrm{L}_{\hat{\nabla}}$ in the second factor is 
$\hat{\nabla}^* \cdot \hat{\nabla}$, which is a self-adjoint elliptic operator. Since $\mathrm{Ker}~\hat{\nabla}=0$, then it follows by standard elliptic theory $\hat{\nabla}^* \cdot \hat{\nabla}$ is invertible on $S^p_{\pm}$. 
By the implicit function theorem in Hilbert space, we obtain a neibourhood 
$U(\omega,\delta)$ with $0<\delta \leq \epsilon$ and a map $\psi:~U(\omega,\delta) \rightarrow P\Gamma(\Xi_{\sigma})$ satisfying our conditions. Hence we 
obtain the diffeomorphism $\Phi$ on $\nabla+U(0,\delta)$ cited in this lemma. 
Since we see that
\[\mathrm{L}_{\hat{\nabla}}([v]^*\omega,[v^{-1}u])=\mathrm{L}_{\hat{\nabla}}(\omega,[u])~,\]
\noindent
for all $[u],[v] \in P\Gamma(\Xi_{\sigma}),~\omega \in \Omega^1(ad_{\sigma}(N_{\tilde{\theta}}))$, then $\Phi$ is chosen $P\Gamma(\Xi_{\sigma})$-equivariantly.
~~~Q.E.D.  \\

\noindent
By $(2)$ of the above lemma 3.18, $U(0,\delta)$ can be chosen as $P\Gamma(\Xi_{\sigma})$-invariant. We then prove the following lemma: \\

\large{\bf{Lemma 3.19}}~~Let $\nabla \in \mathrm{C}_{\mp}(\Xi_{\sigma})$. Then 
there exists a $\delta >0$ such that 
\[\nabla_{\sigma}+\{U(0,\delta)\cap \Omega^1_{\pm}(ad_{\sigma}(N_{\tilde{\theta}}))) \cap \mathrm{Ker}~(\hat{\nabla})^*\} \]
\noindent
is diffeomorphically imbedded in the moduli space 
$\mathcal{M}_{\pm}(\Xi_{\sigma})$ of $\mathcal{C}_{\pm}(\Xi_{\sigma})$ by $\Gamma(\Xi_{\sigma})$ respectively. \\

Proof.~~By Lemma 3.18, $\nabla_{\sigma}+\{U(0,\delta)\cap \Omega^1(ad_{\sigma}(N_{\tilde{\theta}})\}$ is $P\Gamma(\Xi_{\sigma})$-equivariantly diffeomorphic to  $U(0,\mathrm{Id}) \subseteq \mathrm{Ker}~\hat{\nabla}^* \times P\Gamma(\Xi_{\sigma})$. Then we conclude that 
\[(\nabla_{\sigma}+\{U(0,\delta) \cap \mathrm{Ker}~(\hat{\nabla})^*\}) \cap \mathrm{C}_{\pm}(\Xi_{\sigma})\]
\noindent
is diffeomorphically embedded in $\mathcal{M}_{\pm}(\Xi_{\sigma})$. ~~Q.E.D.\\

\noindent
We next show the following lemma (cf:[FU]): \\

\large{\bf{Lemma 3.20}}~~$\mathcal{M}_{\pm}(\Xi_{\sigma})$ is a Hausdorff space.\\

Proof.~~It suffices to show that $\mathcal{M}(\Xi_{\sigma})=\mathcal{C}(\Xi_{\sigma})/\Gamma(\Xi_{\sigma})$ is Hausdorff. So we need to check that 
\[ \Upsilon=\{(\nabla,\gamma_u(\nabla))~|~\nabla \in \mathcal{C}(\Xi_{\sigma}),
u \in \Gamma(\Xi_{\sigma})~\}  \]
\noindent
is closed in $\mathcal{C}(\Xi_{\sigma}) \times \mathcal{C}(\Xi_{\sigma})$. 
In fact, if $\{(\nabla+\omega_j,\gamma_{u_j}(\nabla+\omega_j))\}_j \subseteq \Upsilon$ is convergent to $(\nabla+\omega,\nabla+\omega\prime)$, then 
 $\omega_j \rightarrow \omega$ and 
\[\gamma_{u_j}(\nabla+\omega_j)=\nabla+u_j\hat{\nabla}(u_j^*)+u_j{\omega}u_j^* 
\rightarrow \nabla+\omega\prime ~.\]
\noindent
 Put $\omega\prime_j=u_j\hat{\nabla}(u_j^*)+u_j{\omega}u_j^* \in \Omega^1(ad_{\sigma}(N_{\tilde{\theta}}))$. Then $\hat{\nabla}(u_j^*) =u_j^*\omega\prime-{\omega}u_j^*$. Since $||u_j\nabla||_k~(0\leq k \leq p-1)$ are bouded by definition 
of the Sobolev norms $||\cdot||$, Hence $\{u_j\}$ is bounded in 
$S^p(\Omega^0(ad_{\sigma}(N_{\tilde{\theta}})))$. Then it follows by Rellich's 
theorem that there exists a subsequence $\{u_k\}$ of $\{u_j\}$ such that 
$\{u_k\}$ is convergent in $S^{p-1}(\Omega^1(ad_{\sigma}(N_{\tilde{\theta}}))$. Hence $u_k^*\omega\prime-{\omega}u_k^*$ converges in $S^{p-1}(\Omega^1(ad_{\sigma}(N_{\tilde{\theta}})))$, so is $\{\hat{\nabla}(u_k)\}_k$ in $S^{p-1}(\Omega^1(ad_{\sigma}(N_{\tilde{\theta}})))$ as well. Therefore $\{u_k\}_k$ converges to 
$u \in S^p(\Omega^0(ad_{\sigma}(N_{\tilde{\theta}})))$. Since $\Gamma(\Xi_{\sigma})$ is closed in $S^p(\Omega^0(ad_{\sigma}(N_{\tilde{\theta}})))$, we conclude that $\gamma_u(\nabla+\omega)=\nabla+\omega\prime$. This completes the proof.~~Q.E.D. \\

We now state our main theorem as follows:\\

\large{\bf{Theorem 3.21}}~~Let $M$ be a compact spin 4-dimensional Riemannian 
manifold whose isometry group $\mathrm{Iso}(M,g)$ contains $T^2$, and 
$G \stackrel{\varrho} \rightarrow P \rightarrow M$ a principal 
$G$-bundle over $M$.where $G$ is a compact connected semisimple Lie group. 
Suppose $G \stackrel{\tilde{\varrho}} \rightarrow P_{\tilde{\theta}} \rightarrow M_{\theta}$ exists, 
then given the finite generated projective irreducible right $M_{\theta}$-module $\Xi_{\sigma}=(P_{\theta}\otimes V_{\sigma})^{\tilde{\varrho}\otimes Ad(\sigma)(G)}$ for a highest weight $\sigma$ of $G$,~
the moduli space $\mathcal{M}_{\pm}(\Xi_{\sigma})$ of $\mathcal{C}_{\pm}(\Xi_{\sigma})$ under $\Gamma(\Xi_{\sigma})$ 
is a locally smooth manifold with its dimension:
\[<[M_{\theta}],\mathrm{ch}_{\theta}(\Xi_{\sigma} \otimes S^{\pm})> \pm \mathcal{Rank}_{M_{\theta}}[\mathbb{H}^2_{\mp}] \]
\noindent
respectively, where $\mathrm{ch}_{\theta}$ is the Connes-Chern character from \\ $K_0(M_{\theta})$ to the perodic cyclic homology of $M_{\theta}$, and $S^{\pm}$ are the smooth sections of the $(\pm)$-half spin structure of $M$ respectively.\\
Proof. ~~By~Corollary 3.9,~Lemma 3.17,~3.19 and 3.20,~$\mathcal{M}_{\pm}(\Xi_{\sigma})$ is a smooth manifold.  By Lemma 3.16, its dimension is equal to 
$\mathrm{rank}_{M_{\theta}}\mathrm{H}^1_{\pm}$, which is by Corollary 3.9 and 
Lemma 3.10 that 
\[<[M_{\theta}],\mathrm{ch}_{\theta}(\Xi_{\sigma} \otimes S^{\pm})> \pm \mathcal{Rank}_{M_{\theta}}[\mathbb{H}^2_{\mp}]~.\]
\noindent
respectively. This completes the proof. ~~~Q.E.D. \\

\noindent
From this theorem, we deduce several useful corollaries which are well known in  the case of ordinary manifolds as well as a noncommutative 4-sphere case:\\

\large{\bf{Corollary 3.22}}~~In the above theorem, suppose $[\mathbb{H}^2_{\mp}]=0 \in K_0(M_{\theta})$ for all $\nabla \in \mathcal{C}_{\pm}(\Xi_{\sigma})$, then the instanton moduli space $\mathcal{M}_{\pm}(\Xi_{\sigma})$ of $Xi_{\sigma}$ is a smooth manifold with its dimension:
\[<[M_{\theta}],\mathrm{ch}_{\theta}(\Xi_{\sigma} \otimes S^{\pm})> \]
\noindent
respectively. \\

Proof. ~~By Corollary 3.12 and Theorem 3.21, the conclusion follows. ~~~Q.E.D.\\

\large{\bf{Corollary 3.23}}~~Let $M$ be a compact self-dual(resp.anti self-dual) Riemannian 4-manifold with positive scalar curvature,and G a connected compact semisimple  Lie group. Suppose there exists a noncommutative principal 
$G$-bundle 
\[ G \stackrel{\tilde{\varrho}}\rightarrow N_{\tilde{\theta}} \rightarrow M_{\theta} \]
\noindent
associated with a $T^2$ isometric action on $M$, then given a highest weight 
$\sigma$ of $G$ and consider the right $M_{\theta}$-module 
\[\Xi_{\sigma}=(N_{\tilde{\theta}} \otimes V_{\sigma})^{\tilde{\varrho}
\otimes \sigma^{-1}(G)}  \]~,
\noindent
then $\mathcal{M}_{\pm}(\Xi_{\theta})$ is a smooth manifold with dimension:
\[ \mathrm{dim}\mathcal{M}_{\pm}(\Xi_{\sigma})=
\pm p_1(ad_{\sigma}(N))-\frac{1}{2}\mathrm{dim}G(\chi(M)\mp \tau(M)) \]
\noindent
respectively, where $p_1(ad_{\sigma}(N))$ is the 1-Pontrjagin number of the 
adjoint bundle 
\[ad_{\sigma}(N)=(C^{\infty}(N)\otimes \mathcal{SU}(n_{\sigma})_{\mathbb{C}})^
{\varrho \otimes Ad(\sigma)(G)} \]
\noindent
associated with $\sigma$, $\chi(M)$ is the Euler characteristic of $M$ and 
$\tau(M)$ is the Hirzebruch signature of $M$.\\

Proof.~By the same method as in [AHS], we have that $(\mathbb{H}^0)^2_{\pm}=0$ 
for the 2-cohomology $(\mathbb{H}^0)^2_{\pm}$ with respect to the Laplacian 
associated with $(d_M)^*+P^0_{\mp}\cdot d_M$. Since the quantization map $L_{\theta}$ from $C^{\infty}(M)$ to $M_{\theta}$ can be lifted to a K-theoretic isomorphism (cf:[R]) and $L_{\theta} \cdot d_M = d_{\theta} \cdot L_{\theta}$, then it implies that $\mathbb{H}^2_{\mp}=0$ for all $\nabla \in \mathcal{C}_{\pm}(\Xi_{\sigma})$. Hence it follows from Corollary 3.22 that 
$\mathcal{M}_{\pm}(\Xi_{\sigma})$ is a smooth manifold with their dimension:
\[ \mathrm{dim}\mathcal{M}_{\pm}(\Xi_{\sigma})=
<[M_{\theta}],\mathrm{ch}_{\theta}(\mathrm{P}_{\sigma})\cdot \mathrm{ch}_{\theta}(S^{\pm}_{\theta})> \]
\noindent
Since the cyclic homology (resp.cohomology) of $M_{\theta}$ is isomorphic to 
the de Rham cohomology (resp.homology) of $M$ via the split map $L_{\theta}$(cf:[R]), it follows from [CV] (cf:[C],[MS]) that 
\[<[M_{\theta}],\mathrm{ch}_{\theta}(\mathrm{P}_{\sigma})\cdot \mathrm{ch}_{\theta}(S^{\pm}_{\theta})>=
<[M],\mathrm{ch}(\mathrm{P}^0_{\sigma})\cdot \mathrm{ch}(S^{\pm})>~, \]
\noindent
which is by [AHS] equal to 
\[\pm p_1(ad_{\sigma}(N))-\frac{1}{2}\mathrm{dim}G(\chi(M)\mp \tau(M))~.\]
\noindent
~~~Q.E.D. \\

\large{\bf{Corollary 3.24}}([LS])~~Let us take the Hoph bundle:
\[ \mathrm{SU}(2) \longrightarrow S^7 \longrightarrow S^4~, \]
\noindent
and its associated noncommutative Hoph bundle:
\[ \mathrm{SU}(2) \longrightarrow S^7_{\theta} \longrightarrow S^4_{\theta} \]
\noindent
Let $\Xi_1$ be the right $S^4_{\theta}$-module associated with the highest 
weight 1 of $\mathrm{SU}(2)$ and let $\mathcal{M}_{+}(\Xi_1)$ be the moduli 
space of $\mathcal{C}_{+}(\Xi_1)$. Then it is a smooth manifold with 
its dimension 5.  \\

\large{\bf{Remark}}~~One of the most important problem in this area is to find 
(anti)instantons for a given noncommutaive principal bundle and its associated  noncommutative vector bundles. Even in commutative cases, there exist no (anti)
instanton due to [AHS], in which they presented the complete description of existence of (anti)instantons for principal compact simple linear Lie group bundles over 4-sphere. In the forthcoming paper, we shall discuss this problem for 
noncommutative principal bundles.\\

\begin{center}
\Large{References}
\end{center}
\noindent
[AHS]~M.F.Atiyah,~N.J.Hitchin and I.M.Singer,~Self duality in 4-\\
dimensional Riemannian geometry,~Proc.R.S.London,A362,425-461 (1978). \\

\noindent
[C]~A.Connes,~Noncommutative Geometry,~Academic Press,(1994). \\

\noindent
[CL]~A.Connes and G.Landi,~Noncommutative manifolds: The instanton algebras 
and isospectral deformations,~Comm.Math.Phy.,\\
221,141-159 (2001). \\

\noindent
[CV]~A.Connes and M.Dubois-Violette,~Noncommutative Finite-\\
Dimensional Manifolds.I.~Spherical Manifolds and Related Examples,~
Commun.Math.Phy.230,539-579 (2002). \\

\noindent
[FU]~D.S.Freed and K.K.Uhlenbeck,~Instantons and 4-Manifolds,\\
MSRI.Publ.1,Springer-Verlag (1984).\\

\noindent
[LS]~G.Landi and W.van Suijlekom,~Noncommutative instantons from twisted 
conformal symmetries,~arXiv:math.QA/0601554v2.\\

\noindent
[MS]~C.C.Moore and C.Schochet,~Global Analysis on Foliated Spaces,\\
MSRI.Publ.9,Springer-Verlag (1988).\\

\end{document}